\documentclass[10pt]{article}
\usepackage{amssymb, amsmath, url, graphicx, setspace, geometry}
\usepackage{calrsfs}
\usepackage{wasysym}

\def\3{\subset }
\def\4{\subseteq }
\def\<{\left<}
\def\>{\right>}

\def\bit{\begin{itemize}}
\def\eit{\end{itemize}}
\def\3{\subset }
\def\4{\subseteq }

\def\0{\leqno}

\def\barr{\begin{array}}
\def\earr{\end{array}}

\def\Z{{\rlap{$\kern2pt{\rm Z}$}{\rm Z}\,}}
\def\bld#1#2{{\buildrel{#1}\over{#2}}}
\def\st#1#2{{\mathrel{\mathop{#2}\limits_{#1}}{}\!}}
\def\stb#1#2#3{{\st{{#1}}{\bld{{#2}}{#3}}{}\!}}
\def\xmare#1#2{\stb{#1}{#2}{\mbox{\Huge$\times$}}}

\title{\bf A density result on the sum of element orders of a finite group}
\author{Mihai-Silviu Lazorec and Marius T\u arn\u auceanu}
\date{November 24, 2019}

\begin{document}

\maketitle

\begin{abstract}
Let $\mathcal{G}$ be the class of all finite groups and consider the function $\psi'':\mathcal{G}\longrightarrow(0,1]$, given by $\psi''(G)=\frac{\psi(G)}{|G|^2}$, where $\psi(G)$ is the sum of element orders of a finite group $G$. In this paper, we show that the image of $\psi''$ is a dense set in $[0, 1]$. Also, we study the injectivity and the surjectivity of $\psi''$.
\end{abstract}



\section{Introduction}

One of the main research directions in what concerns the finite group theory is to introduce new tools that are used to characterize the nature (cyclic, abelian, nilpotent, supersolvable, solvable) of a finite group $G$. For an element $x$ of $G$, we denote the order of $x$ by $o(x)$. Approximately one decade ago, in \cite{1}, Amiri, Jafarian Amiri and Isaacs initiated the study of the quantity 
$$\psi(G)=\sum\limits_{x\in G}o(x),$$
which is the sum of element orders of a finite group $G$. Their main result states that, given a group $G$ of order $n$, we have $\psi(G)\leq \psi(C_n)$, and the equality holds if and only if $G\cong C_n$, where $C_n$ is the cyclic group of order $n$.

The sum of element orders was used to provide some criteria which indicate that once the ratio $\psi'(G)=\frac{\psi(G)}{\psi(C_{|G|})}$ is sufficiently large, then $G$ belongs to one of the well-known classes of finite groups. For instance, in \cite{6, 8}, Herzog, Longobardi and Maj proved that if $\psi'(G)>\frac{7}{11}$, then $G$ is cyclic. In \cite{13}, T\u arn\u auceanu showed that once $\psi'(G)>\frac{13}{21}$, $G$ is a nilpotent group. The author conjectures that if $\psi'(G)>\frac{31}{77}$, then $G$ is supersolvable, this result being completely proved by Baniasad Azad and Khosravi in \cite{5}. Finally, in \cite{7}, Herzog, Longobardi and Maj established that once $\psi'(G)>\frac{1}{6.68}$, $G$ is solvable. They outlined the fact that their result may be further improved by replacing the lower bound $\frac{1}{6.68}$ with $\frac{211}{1617}$. Indeed, in \cite{4}, Baniasad Azad and Khosravi showed that $\psi'(G)>\frac{211}{1617}$ implies the solvability of $G$. Note that the lower bounds $\frac{7}{11}$, $\frac{13}{21}$, $\frac{31}{77}$ and $\frac{211}{1617}$ are the best possible ones and the equality between $\psi'(G)$ and each of the previous ratios was also studied in the above mentioned papers. 

In \cite{12}, T\u arn\u auceanu considered the ratio $$\psi''(G)=\frac{\psi(G)}{|G|^2},$$ which is a quantity that is clearly related to the sum of element orders of a finite group $G$. Using this ratio, in the same paper, the author provides some criteria that imply the cyclicity, the abelianess, the nilpotency, the supersolvability and the solvability of $G$ via $\psi''(G)$.

Let $\mathcal{G}$ be the class of all finite groups. It is obvious that $\psi''(G)\in (0,1]$, for all $G\in\mathcal{G}$. Therefore we consider the function $$\psi'':\mathcal{G}\longrightarrow (0, 1]$$ and we study some of its properties. 

Our main objective is to show that the set containing the ratios $\psi''(G)$, where $G\in\mathcal{G}$, is dense in $[0, 1].$ In other words, we prove that the image (range) of $\psi''$ is dense in $[0, 1]$, i.e.:\\

\textbf{Theorem 1.1.} \textit{The set $$Im \ \psi''=\lbrace \psi''(G) \ | \ G\in\mathcal{G}\rbrace$$ is dense in $[0,1]$.}\\
   
We denote by $\psi''|_{\mathcal{S}}$ the restriction of $\psi''$ to a subclass $\mathcal{S}$ of $\mathcal{G}$ containing all finite groups that satisfy a certain property. If we denote by $\mathcal{C}$ the class of all finite cyclic groups, we can write the following consequence of the proof of Theorem 1.1.\\

\textbf{Corollary 1.2.} \textit{The set $Im \ \psi''|_{\mathcal{C}}$ is dense in $[0, 1].$}\\

Further, if we let $\mathcal{S}$ to be a class of finite groups such that $\mathcal{C}\subsetneq \mathcal{S}\subsetneq\mathcal{G}$, we have
$$Im \ \psi''|_{\mathcal{C}}\subseteq Im \ \psi''|_{\mathcal{S}}\subseteq Im \ \psi''.$$
Considering the closures of these sets, Corollary 1.2 and Theorem 1.1 imply the following result.\\

\textbf{Corollary 1.3.} \textit{Let $S$ be a class of finite groups such that $\mathcal{C}\subsetneq \mathcal{S}\subsetneq \mathcal{G}$. Then the set $Im \ \psi''|_{\mathcal{S}}$ is dense in $[0, 1].$}\\

Further, we study other properties of $\psi''$ and we focus especially on the surjectivity and the injectivity of this function. Firstly, we recall the following two results related to the sum of element orders:
\begin{itemize}
\item[--] $\psi(C_{p^{\alpha}})=\frac{p^{2\alpha+1}+1}{p+1}$, where $p$ is a prime and $\alpha\geq 1$ is an integer (see Lemma 2.9 (1) of \cite{6});
\item[--] $\psi$ is multiplicative, i.e. if $G_1$ and $G_2$ are two finite groups such that $(|G_1|, |G_2|)=1$ then $\psi(G_1\times G_2)=\psi(G_1)\psi(G_2)$  (see Lemma 2.1 of \cite{2}).
\end{itemize}

The first result implies that
$$\psi''(C_{p^{\alpha}})=\frac{p^{2\alpha+1}+1}{p^{2\alpha}(p+1)},$$
while the second one guarantees the multiplicativity of $\psi''$. Hence, we have
$$\psi''(C_n)=\prod\limits_{i=1}^k\psi''(C_{p_i^{\alpha_i}})=\prod\limits_{i=1}^k\frac{p_i^{2\alpha_i+1}+1}{p_i^{2\alpha_i}(p_i+1)},$$
where $n=p_1^{\alpha_1} p_2^{\alpha_2}\ldots p_k^{\alpha_k}$ is the prime factorization of a positive integer $n\geq 2$. It is clear that $$\frac{p_i}{p_i+1}<\psi''(C_{p_i^{\alpha_i}})\leq \frac{p_i^3+1}{p_i^3+p_i^2}, \ \forall \ i\in\lbrace 1, 2,\ldots, k\rbrace,$$
the second inequality being a consequence of the fact that the function $f:[1, \infty)\longrightarrow (0, 1)$, given by $f(x)=\frac{p^{2x+1}+1}{p^{2x}(p+1)}$, is strictly decreasing and $f(1)=\frac{p^3+1}{p^3+p^2}$, for any prime $p$. It follows that
$$\prod\limits_{i=1}^k\frac{p_i}{p_i+1}<\psi''(C_n)\leq\prod\limits_{i=1}^k\frac{p_i^3+1}{p_i^3+p_i^2},$$
and, moreover, we have 
$$\psi''(C_n)=\prod\limits_{i=1}^k\frac{p_i^3+1}{p_i^3+p_i^2}\Longleftrightarrow \alpha_i=1, \ \forall \ i\in\lbrace 1, 2, \ldots, k\rbrace \ \text{(i.e. $C_n$ is of square-free order)}.$$

Let $a$ and $r$ be some positive odd integers. We remark that the ratios $\frac{a}{2^r}$ which are greater than $\frac{7}{16}$ are not contained in $Im \ \psi''$. Indeed, if there is a finite group $G$ such that $\psi''(G)=\frac{a}{2^r}$, since $\psi''(G)>\frac{7}{16}$, using Theorem 1.1 c) of \cite{12}, we deduce that $G$ is cyclic. If $G$ is trivial, we obtain $\frac{a}{2^r}=1$, which leads to a contradiction since $a$ is odd. If $G$ is not trivial, there is a positive integer $n\geq 2$ such that $G\cong C_n$. Let $n=p_1^{\alpha_1} p_2^{\alpha_2}\ldots p_k^{\alpha_k}$ the prime factorization of $n$, where $p_1<p_2<\ldots<p_k$. Then, we have
$$\psi''(G)=\frac{a}{2^r}\Longleftrightarrow\prod\limits_{i=1}^k\frac{p_i^{2\alpha_i+1}+1}{p_i^{2\alpha_i}(p_i+1)}=\frac{a}{2^r}\Longleftrightarrow 2^r\prod\limits_{i=1}^k\frac{p_i^{2\alpha_i+1}+1}{p_i+1}=a\prod\limits_{i=1}^k p_i^{2\alpha_i}.$$
Since $a$ is odd, it follows that 2 divides the product $\prod\limits_{i=1}^k p_i^{2\alpha_i}$, so $p_1=2$. Hence, we have
$$2^r\prod\limits_{i=1}^k\frac{p_i^{2\alpha_i+1}+1}{p_i+1}=2^{2\alpha_1}a\prod\limits_{i=2}^kp_i^{2\alpha_i}.$$    
The products $\prod\limits_{i=1}^k\frac{p_i^{2\alpha_i+1}+1}{p_i+1}$ and $a\prod\limits_{i=2}^kp_i^{2\alpha_i}$ are positive odd integers and this leads us to $2^{r}|2^{2\alpha_1}$ and $2^{2\alpha_1}|2^{r}$. Consequently, $r=2\alpha_1$, contradicting the fact that $r$ is odd. This paragraph justifies the following result.\\

\textbf{Proposition 1.4.} \textit{$\psi''$ is not a surjective function.}\\

Let 
$$E(27)=\langle x,y \ | \ x^3=y^3=[x,y]^3=1, [x,y]\in Z(E(27))\rangle,$$
be the non-abelian group of order 27 and exponent 3. Since $$\psi''(C_3^3)=\psi''(E(27))=\frac{79}{729}$$
and $C_3^3\not\cong E(27)$, the next result easily follows.\\

\textbf{Proposition 1.5.} \textit{$\psi''$ is not an injective function.}\\

However, for an arbitrary finite group $G$, there are finite groups $G_0$ such that $\psi''(G)=\psi''(G_0)$ implies $G\cong G_0$. For instance, if $\psi''(G)=\psi''(C_2)$, then $G\cong C_2$. Indeed, since $\psi''(G)=\frac{3}{4}$ and $\frac{3}{4}>\frac{7}{16}$, there is a positive integer $n\geq 2$ such that $G\cong C_n$. Let $n=p_1^{\alpha_1} p_2^{\alpha_2}\ldots p_k^{\alpha_k}$ the prime factorization of $n$, where $p_1<p_2<\ldots<p_k$. One may repeat one of our previous reasonings to show that $p_1=2$. Hence, we have
$$\frac{3}{4}=\psi''(G)=\psi''(C_{2^{\alpha_i}})\prod\limits_{i=2}^k\psi''(C_{{p_i}^{\alpha_i}})\leq\frac{3}{4}\prod\limits_{i=2}^k\psi''(C_{{p_i}^{\alpha_i}})\leq\frac{3}{4},$$
which leads us to $\alpha_1=1$ and $\alpha_i=0$, for all $i\in\lbrace 2, 3, \ldots, k\rbrace$. Therefore, $n=2$, so $G\cong C_2$.

We end this section by outlining an open problem concerning the injectivity of $\psi''$ once we work only with finite cyclic groups.\\

\textbf{Open problem.} \textit{Study the injectivity of the function $\psi''|_{\mathcal{C}}$.}\\

We mention that we wrote a code in Microsoft Visual Studio 2019 \cite{10} and, up to $10^6$, we were not able to find two distinct positive integers $m$ and $n$ such that $\psi''(C_m)=\psi''(C_n)$. Hence, we would say that $\psi''|_{\mathcal{C}}$ is injective.
  
\section{Proof of the main result}

In this section, we prove the validity of Theorem 1.1. First of all, we recall the following preliminary result which is a consequence of the Proposition outlined on the page 863 of \cite{11}.\\

\textbf{Lemma 2.1.} \textit{Let $(x_n)_{n\geq 1}$ be a sequence of positive real numbers such that $\displaystyle \lim_{n \to\infty}x_n=0$ and $\sum\limits_{n=1}^{\infty}x_n$ is divergent. Then the set containing the sums of all finite subsequences of $(x_n)_{n\geq 1}$ is dense in $[0,\infty)$.}\\

A proof of this result is given in \cite{9} (see Lemma 4.1), where Lazorec also studies the density of a set formed of some algebraic quantities. Also, we recall that the series $\sum\limits_{i=1}^{\infty}\frac{1}{p_i}$ is divergent, where $p_i$ is the $i$th prime number. In addition, it is known that a continuous function $f:(X, \tau_X)\longrightarrow (Y, \tau_Y)$, where $(X, \tau_X)$ and $(Y, \tau_Y)$ are topological spaces, has the following property: for any two subsets $A$ and $B$ of $X$, if $\overline{A}_{\tau_X}=\overline{B}_{\tau_X}$, then $\overline{f(A)}_{\tau_Y}=\overline{f(B)}_{\tau_Y}$. Here, we denoted by $\overline{S}_{\tau}$ the closure of a set $S$ with respect to a topology $\tau$. However, when we work with the usual topology $\tau_{\mathbb{R}}$ of $\mathbb{R}$, instead of writing $\overline{S}_{\tau_{\mathbb{R}}}$, we simply write $\overline{S}$. Finally, we outline the following well-known criterion related to the nature of two series (see Theorem 10.9 of \cite{3}).\\ 

\textbf{Lemma 2.2.} \textit{Let $(a_n)_{n\geq 1}$ and $(b_n)_{n\geq 1}$ be two sequences of positive real numbers such that $\displaystyle\lim_{n\to\infty}\frac{a_n}{b_n}=\alpha\in [0, \infty]$. If $\alpha\in (0, \infty)$, then the series $\sum\limits_{n=1}^{\infty}a_n$ and $\sum\limits_{n=1}^{\infty}b_n$ have the same nature.}\\

\textbf{Proof of Theorem 1.1.} Consider the sequence $(k_n)_{n\geq 1}\subset Im \ \psi''|_{\mathcal{C}}$, where $$k_n=\psi''|_{\mathcal{C}}\bigg(\xmare{i\in I}{ }C_{p_i^n}\bigg)=\psi''\bigg(\xmare{i\in I}{ }C_{p_i^n}\bigg), \text{$I$ is a finite subset of $\mathbb{N}^*$ and $p_i$ is the $i$th prime number}.$$ Then,
$$\displaystyle\lim_{n\to\infty}k_n=\displaystyle\lim_{n\to\infty}\psi''\bigg(\xmare{i\in I}{ }C_{p_i^n}\bigg)=\displaystyle\lim_{n\to\infty}\prod\limits_{i\in I}\frac{p_i^{2n+1}+1}{p_i^{2n}(p_i+1)}=\prod\limits_{i\in I}\frac{p_i}{p_i+1}.$$
Hence, we have 
$$\bigg\lbrace \prod\limits_{i\in I}\frac{p_i}{p_i+1} \ \bigg| \ I\subset \mathbb{N}^*, |I|<\infty, \ p_i=\text{$i$th prime number}\bigg\rbrace\subseteq \overline{Im \ \psi''|_{\mathcal{C}}}\subseteq \overline{Im \ \psi''} \subseteq [0, 1].$$
Consequently, if we show that the first of the above sets is dense in $[0, 1]$, then Theorem 1.1 and Corollary 1.2 both hold. Hence, in what follows, we prove that 
$$\overline{\bigg\lbrace \prod\limits_{i\in I}\frac{p_i}{p_i+1} \ \bigg| \ I\subset \mathbb{N}^*, |I|<\infty, \ p_i=\text{$i$th prime number}\bigg\rbrace}=[0, 1].$$

Consider the sequence $(x_i)_{i\geq 1}\subset (0, \infty)$, where $x_i=\ln(\frac{p_i+1}{p_i})$. We have 
$$\displaystyle\lim_{i\to\infty}\frac{x_i}{\frac{1}{p_i}}=\displaystyle\lim_{i\to\infty}\frac{\ln(\frac{p_i+1}{p_i})}{\frac{1}{p_i}}=1.$$
Therefore, since the series $\sum\limits_{i=1}^{\infty}\frac{1}{p_i}$ is divergent, by Lemma 2.2, we deduce that the series $\sum\limits_{i=1}^{\infty}x_i$ is also divergent. It is obvious that $\displaystyle\lim_{i\to\infty}x_i=0$, so all hypotheses of Lemma 2.1 are satisfied. Therefore, we have
$$\overline{\bigg\lbrace \sum\limits_{i\in I}x_i \ \bigg| \ I\subset \mathbb{N}^*, |I|<\infty\bigg\rbrace}=[0, \infty)\Longleftrightarrow$$ $$ \overline{\bigg\lbrace \sum\limits_{i\in I}\ln\bigg(\frac{p_i+1}{p_i}\bigg) \ \bigg| \ I\subset \mathbb{N}^*, |I|<\infty, p_i=\text{$i$th prime number}\bigg\rbrace}=[0, \infty)\Longleftrightarrow$$
$$\overline{\bigg\lbrace \ln\bigg(\prod\limits_{i\in I}\frac{p_i+1}{p_i}\bigg) \ \bigg| \ I\subset \mathbb{N}^*, |I|<\infty, p_i=\text{$i$th prime number}\bigg\rbrace}=[0, \infty).$$

Further, we denote the interval $(0, \infty)$ by $Y$. Consider the topological spaces $(\mathbb{R}, \tau_{\mathbb{R}})$ and $(Y, \tau_Y)$, where $\tau_{\mathbb{R}}$ is the usual topology of $\mathbb{R}$ and $\tau_Y$ is the subspace topology on $Y$. Note that for a subset $S$ of $\mathbb{R}$, we have $\overline{S}_{\tau_Y}=\overline{S}\cap Y$. Since the function 
$$exp:(\mathbb{R}, \tau_{\mathbb{R}})\longrightarrow (\mathbb{R}, \tau_{\mathbb{R}}), \text{ given by } exp(x)=e^x, \ \forall \ x\in\mathbb{R},$$
is continuous and $[1, \infty)$ is a closed set of $\mathbb{R}$, we have
$$\overline{\bigg\lbrace \prod\limits_{i\in I}\frac{p_i+1}{p_i} \ \bigg| \ I\subset \mathbb{N}^*, |I|<\infty, p_i=\text{$i$th prime number}\bigg\rbrace}=[1, \infty).$$
Note that
$$\overline{\bigg\lbrace \prod\limits_{i\in I}\frac{p_i+1}{p_i} \ \bigg| \ I\subset \mathbb{N}^*, |I|<\infty, p_i=\text{$i$th prime number}\bigg\rbrace}_{\tau_Y}=$$
$$\overline{\bigg\lbrace \prod\limits_{i\in I}\frac{p_i+1}{p_i} \ \bigg| \ I\subset \mathbb{N}^*, |I|<\infty, p_i=\text{$i$th prime number}\bigg\rbrace}\cap Y=$$
$$[1,\infty)\cap Y=\overline{[1, \infty)}\cap Y=\overline{[1,\infty)}_{\tau_Y}.$$
Hence, if we consider the continuous function
$$f:Y\longrightarrow \mathbb{R}, \text{ given by } f(y)=\frac{1}{y}, \ \forall \ y\in Y,$$
we deduce that
$$\overline{\bigg\lbrace \prod\limits_{i\in I}\frac{p_i}{p_i+1} \ \bigg| \ I\subset \mathbb{N}^*, |I|<\infty, p_i=\text{$i$th prime number}\bigg\rbrace}=\overline{(0, 1]}=[0, 1].$$
Consequently, our proof is complete.
\hfill\rule{1,5mm}{1,5mm}\\

\vspace*{3ex}
\small

\begin{minipage}[t]{7cm}
Mihai-Silviu Lazorec \\
Faculty of  Mathematics \\
"Al.I. Cuza" University \\
Ia\c si, Romania \\
e-mail: {\tt mihai.lazorec@student.uaic.ro}
\end{minipage}
\hspace{3cm}
\begin{minipage}[t]{7cm}
Marius T\u arn\u auceanu \\
Faculty of  Mathematics \\
"Al.I. Cuza" University \\
Ia\c si, Romania \\
e-mail: {\tt tarnauc@uaic.ro}
\end{minipage}
\end{document}